\documentclass[a4paper,12pt,reqno]{amsart}

\usepackage{amsthm}
\usepackage{amsmath}
\usepackage{amssymb}
\usepackage{mathrsfs}
\usepackage{latexsym}
\usepackage{exscale}
\usepackage{geometry}
\usepackage{graphicx}
\usepackage[utf8]{inputenc}
\usepackage{verbatim}

\usepackage{enumerate}

\usepackage{color}
\definecolor{red}{rgb}{1,0.1,0.1}
\definecolor{blue}{rgb}{0.1,0.1,1}
\definecolor{vb}{RGB}{160,32,240}

\theoremstyle{plain}

\newtheorem*{teo*}{Theorem}
\newtheorem*{prop*}{Proposition}

\numberwithin{equation}{section}
\newtheorem{teo}{Theorem}[section]

\theoremstyle{remark}

\theoremstyle{definition}

\newtheorem*{mydef*}{Definition}

\calclayout

\allowdisplaybreaks

\begin{document}
	\title[harmonic analysis in nilmanifols]{A generalized Gelfand pair attached to a 3-step nilpotent Lie group}
	
	\author[A.~L.~Gallo]{Andrea L. Gallo}
	\address{A.~L.~Gallo \\ FaMAF \\ Universidad Nacional de C\'ordoba \\
		CIEM (CONICET) \\ 5000 C\'ordoba, Argentina}
	\email{andregallo88@gmail.com}

	\author[L.~V.~Saal]{Linda V. Saal}
	\address{L.~~Saal\\ FaMAF \\ Universidad Nacional de C\'ordoba \\
		CIEM (CONICET) \\ 5000 C\'ordoba, Argentina}
	\email{saal@mate.uncor.edu}

	\thanks{ The authors are  partially supported by
		CONICET and SECYT-UNC}

	
	\dedicatory{\today}
	
	

	\begin{abstract}
		Let $N$ be a nilpotent Lie group and $K$ a compact subgroup of the automorphism group $Aut(N)$ of $N$. It is  well-known that if $(K\ltimes N,N)$ is a Gelfand pair then $N$ is at most $2$-step nilpotent Lie group. 
		
		The notion of Gelfand pair was generalized when $K$ is a non-compact group. In this work, we give an example of a $3$-step
		nilpotent Lie group and a non-compact subgroup $K$ of $Aut(N)$ such
		that $(K \ltimes N,N)$ is a generalized Gelfand pair.
	\end{abstract}
	
	\maketitle
	
\section{Introduction}
Let $G$ be a Lie group and $K$ a compact subgroup of $G$. We denote by $\mathcal{D}(G/K) $ the space of $C^{\infty}$-functions on $G/K$
with compact support and by $\mathcal{D}_{K}(G) $ the subspace
of $\mathcal{D}(G) $ of functions on $G$ which are right
$K$-invariant. Both spaces are identified by mapping $f\in \mathcal{D}(G/K) $ to $f_{0}:=f\circ \mathcal{P}$, where $\mathcal{P}:G\rightarrow
G/K$ is the natural projection.

It follows from the Schwartz's kernel Theorem, that every linear operator which maps continously
 $\mathcal{D}(G/K) $ in $\mathcal{D}'(G/K) $ with respect to the standard topologies
and commuting with the action of $G$ is a convolution operator with a $K$-bi-invariant distribution in $\mathcal{D}'(G)$.

In particular, we consider the subalgebra of convolution operators which
kernels are $K$-bi-invariant integrable functions on $G$. When this algebra
is commutative, we can expect a kind of simultaneous  ``diagonalization" of
all these operators. This motivated, in part, the study of Gelfand pairs and
the corresponding spherical analysis. 
In this sense, we begin by introducing the concept of Gelfand pair.
The following statements are equivalent:
\begin{enumerate}[i)]
\item The convolution algebra of $K$-bi-invariant integrable functions on $G$
is commutative.

\item For any irreducible unitary representation $(\pi,\mathcal{H}) $ of $G,$ the subspace $\mathcal{H}_{K}$ of vectors fixed by $K$
is at most one dimensional.
\end{enumerate}

When any of the above holds, we say that $(G,K)$ is a Gelfand pair.\\

Very well studied examples of Gelfand pairs are provided by symetric pairs of compact or non-compact types. More recent works have put attention on Gelfand pairs of the form $(K \ltimes N, K)$ where $N$ is a nilpotent Lie group and $K$ is a subgroup of the automorphism group $Aut(N)$ of $N$ (see \cite{ADF},\cite{BJR1},\cite{BJR2},\cite{BJR3},\cite{FRY},\cite{L},\cite{V}, among others). One of the first results, proved in \cite{BJR1}, stated that if $(K\ltimes N,K)$ (in short $(K,N)$) is a Gelfand pair then $N$ is abelian or a $2$-step nilpotent group.

The notion of Gelfand pair was extended to the case when $K$ is \textit{non-compact}. 
Observe that, in this case, the space of $K$-invariant integrablefunctions on $G/K$ is trivial. 
Thus, one attempt is to generalize $(ii)$. Seminal papers are due to J. Faraut \cite{F} and E.G. Thomas \cite{T}, and
there is a nice survey in \cite{VD}. 
First of all, we assume that $G$
and $K$ are \textit{unimodular} groups.

Let $(\pi,\mathcal{H})$ be a unitary representation of $G$,
and denote by $\mathcal{H}^{\infty}$ the space of $\mathcal{C}^{\infty}$-vectors, that is, $\mathcal{H}^{\infty}=\{v \in \mathcal{H} : g \mapsto \pi(g) v \in \mathcal{C}^{\infty}(G)\}$. $\mathcal{H}^{\infty }$ is a Fr\'{e}chet space equipped with a natural
Sobolev topology. Let $\mathcal{H}^{-\infty}$ be the antidual of $\mathcal{H}^{\infty}$, with the strong topology (uniform convergence on
bounded sets of $\mathcal{H}^{\infty})$. This yields natural imbeddings
\[
{\mathcal{H}}^{\infty }\subset {\mathcal{H}}\subset {\mathcal{H}}^{-\infty
}. 
\]
We denote by $\pi _{\infty}$ the restriction of $\pi$ to $\mathcal{H}^{\infty}$, and for $g\in G$ define $\pi _{-\infty}(g)$ on $\mathcal{H}^{-\infty}$ by duality: for $\phi \in \mathcal{H}^{-\infty }, v\in \mathcal{H}^{\infty}$
 $$\langle \pi _{-\infty}(g)\phi ,v\rangle
:=\langle \phi ,\pi _{\infty}(g)v\rangle.$$
The elements of $\mathcal{H}^{-\infty}$ are called \textit{distribution vectors}.

We say that $(G,K)$ is a \textit{generalized Gelfand pair if
	for any irreducible unitary representation} $(\pi,\mathcal{H})$ of $G$ \textit{the space }$\mathcal{H}_{K}^{-\infty}$ \textit{of
	distribution vectors fixed by} $K$ is at most one dimensional.\bigskip

In this work, we give an example of a $3$-step
nilpotent Lie group and a non-compact subgroup $K$ of $Aut(N)$ such
that $(K,N)$ is a generalized Gelfand pair. This is stated in Theorem \ref{Teo}.\\
 
 \textbf{Acknowledgements:} We are grateful to G. Ratcliff who let us know the article \cite{R}.

\section{Preliminaries}

We begin this section by recalling some known results about generalized Gelfand pairs.

When $(G,K)$ is a Gelfand pair, there is a one-to-one
correspondence between $K$-bi-invariant functions on $G$ of positive type
and equivalent classes of unitary representations having a cyclic vector
fixed by $K$. Moreover, for a $K$-bi-invariant function $\varphi$ on $G$ of
positive type it holds a Bochner-Godement Theorem
\[
\varphi =\int\limits_{\Sigma }\varphi _{\lambda }d\mu \left( \lambda
\right) 
\]
where $\Sigma$ denotes the set of \textit{extremal} $K$-bi-invariant functions
on $G$ of positive type (or \textit{spherical functions of positive type}),
and $d\mu$ is a Radon measure on $\Sigma$. This allows to see that any 
\textit{spherical} representation of $G$ decomposes multiplicity free.

\bigskip

There is an analogous result for a generalized Gelfand pair:
Let $(\pi,\mathcal{H})$ be a unitary representation of $G$
having a distribution vector $\phi \in \mathcal{H}_{K}^{-\infty }$. Then
for $f\in \mathcal{D}(G)$, it is easy to see that $\pi_{\infty}(f) \phi \in \mathcal{H}^{\infty}$ and $T_{\phi}$
defined by 
$$T_{\phi }(f) =\langle \phi,\pi_{\infty}(f)\phi \rangle,$$
is a $K$-bi-invariant distribution of positive type. $T$ is called the 
\textit{reproducing kernel} of $\mathcal{H}$.

Conversely, let $T$ be a $K$-bi-invariant distribution of positive type.
For $f\in $ $\mathcal{D}(G/K),$ let $f_{0}=f\circ \mathcal{P}$. On $\mathcal{D}(G/K)$ let us
consider the scalar product $\langle f,g\rangle =$ $T_{\phi}(f_{0}^{\ast}\ast g_{0})$, where $f_{0}^{\ast}(x) =\overline{f_{0}(x^{-1})}$, and denote by $\mathcal{N}$ the subspace of vectors of lenght zero. The Hilbert subspace $\mathcal{H}$
of $\mathcal{D}'(G/K) $ associated to $T$ is the
completion of $\mathcal{D}(G/K)/\mathcal{N}$, and an
easy computation (using that a Hilbert space $\mathcal{H}$ is identified
with its dual) shows that if $J^{\ast}:\mathcal{D}(G/K)
\rightarrow \mathcal{D}(G/K)/\mathcal{N}$ is the
natural projection and $J:\mathcal{H}\rightarrow\mathcal{D}(G/K)'$ is the dual map, then for $f\in \mathcal{D}(G/K)$, $(J\circ J^{\ast}) f=f_{0}\ast T$.

Thus, we have the following result (for a
detailed proof see \cite{VD}). 

\begin{teo*}[A]
There is a one-to-one correspondence between unitary representations of $G$
having a cyclic distribution vector fixed by $K$ and $K$-bi-invariant
distributions of positive type in $\mathcal{D}'(G)$ (the corresponding representation is realized as an invariant Hilbert
subspace of $\mathcal{D}'(G/K)$).
\end{teo*}

Also the following results hold:
\begin{itemize}
\item Bochner-Godement's theorem: For every $K$-bi-invariant
distribution $T$ of positive type there exists a Radon measure on the set $\Sigma$ of extremal $K$-bi-invariant distributions of positive type, such
that $T=\int_{\Sigma}T_{s}d\mu (s)$;

\item  Every $G$-invariant Hilbert subspace of $\mathcal{D}'(G/K)$ decomposes multiplicity free. 
\end{itemize}

Now, let us consider a unimodular Lie group $H$ such that for any $(\gamma,\mathcal{V}) \in \widehat{H}$, $\gamma(f)$ is a class trace operator for all $f\in \mathcal{D}(H)$ 
(this property holds for a wide class of Lie groups such as
nilpotent or semisimple Lie groups).

Let us consider the pair $(G,K)$ where $G=H\times H$ and $K=diag(H\times H)$ which is naturally identified with $H$. Also 
$G/K$ can be identified with $H$. Let us denote by $(\gamma^{\ast},\mathcal{V}^{\ast})$ the contragradient representation of $(\gamma,\mathcal{V}_{\gamma})$. On the first hand, $\mathcal{V}^{\ast}\otimes \mathcal{V}$ is canonically isomorphic to the Hilbert
subspace $\mathcal{H}_{\gamma }$ of $\mathcal{D}'(H)$ of distributions of the form $f\ast \chi _{\gamma},f\in \mathcal{D}(H)$. On the other hand, $\gamma^{\ast}\otimes\gamma $ corresponds to the representation of $H\times H$ on $\mathcal{D}'(H)$ given by $(h_{1},h_{2})\mapsto L(h_{1}) R(h_{2})$. Thus, $\chi_{\gamma}$ is the reproducing kernel of $\mathcal{H}_{\gamma}$ and clearly $\chi_{\gamma}$ is a distribution vector in $\mathcal{H}_{\gamma}^{-\infty}$ fixed by $H$.

\bigskip

The complete result, due to Mokni and Thomas in \cite{MT} yields
an analogous of a Carcano criterion for Gelfand pair.

\begin{teo*}[B]
Let $(\omega,\mathcal{W}),(\gamma,\mathcal{V})$ be unitary representations of $H$ such that $\gamma$ is
irreducible. Then $\gamma$ appears in the decomposition of $\omega$ into
irreducible components if and only if $\gamma^{\ast}\otimes \omega$ has a distribution vector fixed by $H$ as $(H\times H)$-module.
\end{teo*}

\bigskip 

Let $N$ be a nilpotent Lie group and denote by $\widehat{N}$ the set of equivalent class of irreducible unitary representation of $N$. 
We describe $\widehat{N}$ according to Kirillov's theory. Let $\mathfrak{n}$ be the Lie algebra of $N$. The group $N$ acts on $\mathfrak{n}$ by the
adjoint action $Ad$, and  $N$ acts on $\mathfrak{n}^{\ast}$, the dual space
of $\mathfrak{n}$, by the dual representation $Ad^{\ast}(n)
\Lambda =\Lambda \circ Ad(n^{-1}).$ Fixed a non trivial $\Lambda \in \mathfrak{n}^{\ast},$ let $O_{\Lambda} :=\left\{Ad^{\ast}(n) \Lambda :n\in N\right\} $ be its coadjoint orbit.

From Kirillov's theory it follows that there is a
correspondence between $\widehat{N}$ and the set of coadjoint orbits in $\mathfrak{n}^{\ast}$.
Indeed, let 
\begin{equation} \label{corchetes}
B_{\Lambda}(u,v) :=\Lambda([u,v]), \, u,v\in \mathfrak{n}.
\end{equation}
Let $\mathfrak{m}$ be a maximal isotropic subspace of $\mathfrak{n}$, and
set $M=\exp(\mathfrak{m})$. Defining on $M$ the character $\chi_{\Lambda}(\exp u)=e^{i\Lambda(u)}$, the
irreducible representation corresponding to $O_{\Lambda} $ is
the induced representation $\rho_{\Lambda }:=Ind_{M}^{N}(\chi_{\Lambda})$.

Let $K$ be a subgroup of $Aut(N)$. Given $k \in K$, $\Lambda \in \mathfrak{n}^{\ast}$ there is a new representation of $N$ defined by $\rho_{\Lambda}^{k}(n):=\rho_{\Lambda}(k\cdot n)$. The stabilizer of $\pi_{\Lambda}$ is $K_{\Lambda}:=\{k \ | \ \rho_{\Lambda} \sim \rho_{\Lambda}^{k}\}$. For each  $k \in K_{\Lambda}$, one can choose an intertwining operator $\omega_{\Lambda}(k)$ such that
$\rho_{\Lambda}^{k}(n) = \omega_{\Lambda}(k)\rho_{\Lambda}(n)\omega_{\Lambda}(k^{-1})$ for all  $n \in N$. The map $k \mapsto \omega_{\Lambda}(k)$ is a projective representation of  $K_{\Lambda}$, i.e\@, $\omega_{\Lambda}(k_1 k_2)  = \sigma_{\Lambda}(k_1,k_2) \omega_{\Lambda}(k_1) \omega_{\Lambda}(k_2)$, with $|\sigma_{\Lambda}(k_1,k_2)|=1$ for all $k_1,k_2 \in K_{\Lambda}$. The map $\omega_{\Lambda}$ is called the intertwining representation of $\pi_{\Lambda}$ or metaplectic representation and $\sigma_{\Lambda}$ the multiplier for the projective representation $\omega_{\Lambda}$. 


Here we shall consider a $3$-step nilpotent Lie group introduced in \cite{R} and a certain subgroup $K$ of $Aut(N)$, such that
\begin{enumerate}[(i)]
\item for all $\Lambda \in \mathfrak{n}^{\ast},K_{\Lambda }=K$,

\item $\omega _{\Lambda }$ is a true representation of $K$.
\end{enumerate}

In this situation, Mackey theory asserts that for $\sigma \in \widehat{K}$, $\rho_{\sigma,\Lambda}=\sigma \otimes \omega_{\Lambda}\rho_{\Lambda}$
is an irreducible representation of $K\ltimes N$ and by varying $\sigma \in \widehat{K}$ and $\rho_{\Lambda}\in \widehat{N}$ this construction
exhausts $\widehat{K\ltimes N}$.

It follows from Theorem B that the irreducible representation $\sigma \otimes
\omega_{\Lambda}\rho_{\Lambda}$ has a distribution vector fixed by $K$
if and only if the dual representation $\sigma^{\ast}$ of $K$ appears in
the decomposition into irreducible components of $\omega_{\Lambda}$.

As a consequence, $(K,N)$ will be a generalized Gelfand pair if and
only if $\omega_{\Lambda }$ is multiplicity free.

\bigskip




\section{Example}
The group $N$ to be considered is the simplest case of the family introduced by G. Ratcliff in \cite{R}. Let $H_1$ be the 3-dimensional Heisenberg group with Lie algebra $\mathfrak{h}_1$ whose coordinates are $(x,y,t) \in \mathbb{R}^{3}$
and Lie bracket defined by $$[(x,y,t),(x',y',t')]=(0,0,xy'- yx').$$

 Let $S$ be the subgroup of $Sp(1) \subseteq Aut(H_1)$ consinting of the matrices 
\begin{equation*}
\mathbf{s} =  \left(
\begin{array}{cc}
1 & 0 \\
s & 1   
\end{array}
\right), \ s \in \mathbb{R}.
\end{equation*}
Let us consider the action of $S$ on $H_1$ given by 
$$ \mathbf{s} \cdot (x,y,t) = (x,sx+y,t).$$  
This action gives rise to a semidirect product $N=S \ltimes H_1$ such that
\begin{equation}\label{producto}
(s,x,y,t)(s',x',y',t')=(s+s',x+x',sx'+y+y',t+t'+\frac{y'x+sxx'-x'y}{2}).
\end{equation}
Let $\mathfrak{s}$ be the Lie algebra of $S$. The Lie algebra  $\mathfrak{n}$ associated to $N$ is a 3-step nilpotent Lie algebra with coordinates $(s, x, y, t)$, where $s \in \mathfrak{s}$, $x,y,t \in \mathbb{R}$, and product 
\begin{equation}\label{corchete}
[(s, x, y, t), (s’, x’, y’, t’)] = (0, 0, sx' - s'x, xy' -x' . y);
\end{equation}
its one-dimensional center is $\mathfrak{c}=\{(0,0,0,t) \ | \ t \in \mathbb{R}\}$. 
	
We denote by $Aut_0(N)$ the group of automorphisms of $N$ acting on $\mathfrak{c}$ by the identity. Since the exponential map is the identity  
\begin{equation}\label{aut}
Aut_0(N)=\{k \in  GL(4,\mathbb{R}) : k([u,v])=[k(u),k(v)] \text{ for all } u,v \in \mathfrak{n}, k\mid_{\mathfrak{c}}=I\}.
\end{equation}
Let $\Phi \in Aut_0(N)$, and $\mathcal{B}=\{e_j\}_{j=1}^{4}$ be the cannonical basis of $\mathbb{R}^{4}$. According to \eqref{corchete} and \eqref{aut},   $\Phi$ must satisfy the following relationships
$$\Phi([e_1,e_2])=\Phi(e_3); \ \  \Phi([e_2,e_3])=\Phi(e_4);$$ $$\Phi([e_1,e_3])=\Phi([e_1,e_4])=\Phi([e_2,e_4])=\Phi([e_3,e_4])=\Phi(0)=0.$$
 Thus, we obtain that
\begin{equation*}
Aut_0(N)= \lbrace \left(
\begin{array}{cccc}
r & a & 0 & 0 \\
0 & r^{-\frac{1}{2}} & 0 & 0 \\
d & b & r^{\frac{1}{2}} & 0 \\
e & c & -dr^{-\frac{1}{2}} & 1   
\end{array}
\right) \ : \ r,a,b,c,d,e \in \mathbb{R} \rbrace.
\end{equation*}

We define
\begin{equation*}
A := \lbrace \left(
\begin{array}{cccc}
r & 0 & 0 & 0 \\
0 & r^{-\frac{1}{2}} & 0 & 0 \\
0 & 0 & r^{\frac{1}{2}} & 0 \\
0 & 0 & 0 & 1   
\end{array}
\right) \ : \ r \in \mathbb{R} \rbrace \ \ \text{and}
\end{equation*}
\begin{equation*}
M := \lbrace \left(
\begin{array}{cccc}
1 & a & 0 & 0 \\
0 & 1 & 0 & 0 \\
d & b & 1 & 0 \\
e & c & -d & 1   
\end{array}
\right) \ : \ a,b,c,d,e \in \mathbb{R} \rbrace.
\end{equation*}
Then we have that $A$ is acting by conjugation over $M$ and $Aut_0(N) = A \ltimes M$.\\
Writing the elements in $M$ as $5$-uplas $(a,b,c,d,e)$, we have that the product is given by 
\begin{equation}\label{producto}
(a,b,c,d,e) (a',b',c',d',e') = (a+a',b+b'+da',c+c'+ea'-db',d+d',e+e'-dd')
\end{equation}
Moreover, $M$ is isomorphic to $H \ltimes \mathbb{R}^{3}$, where

 %
\begin{equation*}
H := \lbrace \left(
\begin{array}{cccc}
1 & 0 & 0 & 0 \\
0 & 1 & 0 & 0 \\
d & 0 & 1 & 0 \\
e & 0 & -d & 1   
\end{array}
\right) \ : \ d,e \in \mathbb{R} \rbrace, \text{and}
\end{equation*}
\begin{equation*}
\mathbb{R}^{3}:= \lbrace \left(
\begin{array}{cccc}
1 & a & 0 & 0 \\
0 & 1 & 0 & 0 \\
0 & b & 1 & 0 \\
0 & c & 0 & 1   
\end{array}
\right) \ : \ c,d,e \in \mathbb{R} \rbrace.
\end{equation*}
Indeed, the product on $H$ is $$(0,0,0,d,b)(0,0,0,d',b')=(0,0,0,d+d',e+e'-dd'),$$
and thus we can indentify $H$ with the  subgroup $\{(d,-d,e) \ : \ d,e \in \mathbb{R}\}$ of $H_1$.
By considering the action of $H$ over $\mathbb{R}^{3}$ given by 
$$(d,e) \cdot (a',b',c') = (a',b'+da',c'+ea'-db'),$$ and using \eqref{producto} we obtain that $M=H \ltimes \mathbb{R}^{3}$.\\ 

The subgroup $K$ of $Aut_0(N)$ that we will consider is 
\begin{equation}\label{K}
K = \lbrace k \in \mathbb{R}^{3}: k =\left(
\begin{array}{cccc}
1 & 0 & 0 & 0 \\
0 & 1 & 0 & 0 \\
0 & k_1 & 1 & 0 \\
0 & k_2 & 0 & 1  
\end{array}
\right) \ : \ k_1,k_2 \in \mathbb{R} \rbrace.
\end{equation}

Let $K_1$ (resp. $K_2$) be the subgroup of $K$ whose elements have matrix entry $k_2=0$ (resp. $k_1=0$).

We denote by $(\alpha, \mu, \nu, \lambda)$ the elemnts of $\mathfrak{n}^{\ast}$. The pairing between $\mathfrak{n}$ and $\mathfrak{n}^{\ast}$ is given by $$(\alpha,\mu,\nu,\lambda)[(s,x,y,t)]=\alpha s + \mu x + \nu y + \lambda t.$$
For $\Lambda \in \mathfrak{n}^{\ast}$, it is easy to see that  $K_{\Lambda}=\{ k \in K \ | \ k \cdot \Lambda \in O_{\Lambda}\}$ where $O_{\Lambda}$ is the coadjoint orbit of $\Lambda$. 
Let $X_{\Lambda}  \in \mathfrak{n}$  such that $\Lambda(Y)= \langle Y,X_{\Lambda} \rangle$ for all $Y \in \mathfrak{n}$, hence it follows that $k \cdot \Lambda(Y)=\langle Y,{k^{-1}}^{t}X_{\Lambda} \rangle$.
So
\begin{equation}\label{orbitas}
K_{\Lambda}:=\{k \in K \ | \ k^{t} \cdot X_{\Lambda} \in O_{\Lambda}\},
\end{equation}
where $k^{t}$ denotes the transposed of $k$.

 The \textit{generic orbits}  are those which correspond to the representations with non-zero Plancherel measure and were computed in \cite{R}. They are parametrized by $\Lambda=(\alpha,0,0,\lambda)$ with $\lambda \neq 0$, and if $O_{\alpha,\lambda}$ is  the coadjoint orbit of $\Lambda$, then 
$$O_{\alpha,\lambda}=\{(\alpha-\frac{1}{2\lambda}\nu^{2},\mu, \nu, \lambda) \ | \ \mu, \nu \in \mathbb{R}\}.$$ 
 
 In the case of the non-generic orbits with $\Lambda=(\alpha,\mu,\nu,0)$, by the well-known equality $Ad \circ exp = exp \circ ad$, we obtain that $O_{\Lambda}=\{(\beta,\eta,\nu,0) \ | \ \beta,\eta \in \mathbb{R}\}$. Let  $(0,0,\nu,0)$ be a representative of $O_{(\alpha,\mu,\nu,0)}$, and we denote $O_{(\alpha,\mu,\nu,0)}$ by $O_{\nu}$.

 We now compute explicitly the representation $\rho_{\Lambda}$ corresponding to the orbits $O_{\Lambda}$ for all $\Lambda \in \mathfrak{n}^{\ast}$. We denote $\rho_{\Lambda}$ by $\rho_{\alpha,\lambda}$ in the case $\Lambda=(\alpha,0,0,\lambda)$ and by $\rho_{\nu}$ in the case $\Lambda=(0,0,\nu,0)$.

Fixed $\Lambda = (\alpha,0,0,\lambda)$ with $\lambda \neq 0$, the non-degenerate skew-symetric form associated is
\begin{eqnarray*}
	B_{\Lambda}((s,x,y,t),(s',x',y',t'))&=&(\alpha,0,0,\lambda)([(s,x,y,t),(s',x',y',t')])\\
	&=&(\alpha,0,0,\lambda)(0,0,sx'-s'x,xy'-x'y)\\
	&=&\lambda(xy'-x'y).
\end{eqnarray*}
Thus, a maximal isotropic subespace associated to $\Lambda$ is given by $$\mathfrak{M}_{\Lambda}=\{(s,0,y,t) \in \mathfrak{n} \ | \ s,y,t \in \mathbb{R}\}.$$
The character $\chi_{\Lambda}$ defined on $M_{\Lambda}=\exp(\mathfrak{M}_{\Lambda})$ is $\chi_{\Lambda}(s,0,y,t) = e^{\Lambda(s,0,y,t)} = e^{i(\alpha s+\lambda t)}$, and
$$\rho_{\alpha,\lambda}=Ind_{M_{\Lambda}}^{ N} (\chi_{\Lambda}).$$ 
 We recall that the induced representation is the pair $(\rho_{\alpha,\lambda},H_{\alpha,\lambda})$ where $H_{\alpha,\lambda}$ is the completion of  
 $$\{f \in C_c(N) \ | \ f(nm) = \chi_{\Lambda }(m^{-1}) f(n) \ \text{ for all } m\in M_{\Lambda}, n \in N \},$$ 
 with respect to the inner product $\langle f,g \rangle = \int_{N / M_{\Lambda}} f(u) \overline{g(u)} du$, and the action is given by the regular left translation, that is $(\rho_{\alpha,\lambda}(n)f)(n')=f(n^{-1}n')$, $n,n' \in N$. 
 Notice that setting
 \begin{equation*}
(s,x,y,t) = (0,x,0,0) (s,0,y,t-\frac{xy}{2}),
\end{equation*}
 we can identify $H_{\alpha,\lambda}$ with $L^{2}(\mathbb{R})$ via the map $(0,u,0,0) \mapsto u$. 

Since $$(s,x,y,t)^{-1}=(-s,-x,sx-y,-t),$$ for $f \in H_{\alpha,\lambda}$ we obtain
\begin{eqnarray*}
[\rho_{\alpha,\lambda}(s,0,0,0)f](u)&=&f((s,0,0,0)^{-1}(0,u,0,0))\\
&=&f((-s,0,0,0)(0,u,0,0))\\
&=& f(-s,u,-su,0)\\
&=&f((0,u,0,0)(-s,0,-su,\frac{su^{2}}{2})\\
&=&\chi_{\Lambda}({s,0,su,-\frac{su^{2}}{2}}) f(0,u,0,0)\\
&=& e^{s\alpha - \lambda \frac{su^{2}}{2}} f(u).
\end{eqnarray*}

Analogously, we have $$[\rho_{\alpha,\lambda}(s,0,0,0)f] (u) = e^{i s(\alpha - \frac{\lambda u^{2}}{2})} f(u),$$
$$[\rho_{\alpha,\lambda}(0,x,0,0)f](u) = f(u-x),$$
$$[\rho_{\alpha,\lambda}(0,0,y,0)f] (u) = e^{-i\lambda uy} f(u),$$
$$[\rho_{\alpha,\lambda}(0,0,0,t)f] (u) = e^{i\lambda t} f(u).$$ 
We observe that the representations $\rho_{\alpha,\lambda}$ with $\lambda \neq 0$, are extensions of irreducible representations of $H_1$.\\

We now describe the representations corresponding to non-generic orbits $O_{\nu}$ with $\nu \neq 0$. In this case, $B_{\Lambda}((s,x,y,t),(s',x',y',t'))=\nu(sx'-s'x)$, and a maximal isotropic subspace is again $\mathfrak{M}_{\Lambda}=\{(s,0,y,t) : s,y,t \in \mathbb{R}\}$.
The character associated is $\chi_{\Lambda}(s,0,y,t) = e^{i \nu y}$. With similar computations to the above case, we obtain  
$$[\rho_{\nu}(s,0,0,0)f] (u) = e^{i \nu s u}f(u),$$
$$[\rho_{\nu}(0,x,0,0)f](u) = f(u-x),$$
$$[\rho_{\nu}(0,0,y,0)f] (u) = e^{i \nu y} f(u),$$
$$[\rho_{\nu}(0,0,0,t)f] (u) = f(u).$$ \\

By \eqref{orbitas} we can easily see that  $K_{\Lambda}=K$ for all $\Lambda \in \mathfrak{n}^{\ast}$. Thus, the metaplectic representation $\omega_{\Lambda}$ must satisfy
$$\rho_{\Lambda}^{k}(n) \omega_{\Lambda}(k) = \omega_{\Lambda}(k) \rho_{\Lambda}(n)\quad \text{for all $k \in K$ and $n \in N$}.$$

The subgroups $K_1$ and $K_2$ fix the elements $(s,0,0,0)$, $(0,0,y,0)$ and $(0,0,0,t)$ for all $s,y,t \in \mathbb{R}$. Then, we have to find an unitary operator $\omega_{\Lambda}$ on $L^{2}(\mathbb{R}^{2})$ such that  
\begin{equation}\label{conmut}
\rho_{\Lambda}^{k}(0,x,0,0) \omega_{\Lambda}(k) = \omega_{\Lambda}(k) \rho_{\Lambda}(0,x,0,0), \ \ \forall x \in \mathbb{R}, k \in K
\end{equation}
and 
\begin{equation}\label{otros}
\rho_{\Lambda }(n)\omega_{\Lambda}(k)=\omega_{\Lambda}(k)\rho_{\Lambda}(n)
\end{equation}
 for $n=(s,0,0,0)$, $n=(0,0,y,0)$ and $n=(0,0,0,t)$, $k \in K$.\\
 
 We denote by $\omega_{\alpha,\lambda}$ (resp. $\omega_{\nu}$) the metaplectic representation corresponding to $\Lambda=(\alpha,0,0,\lambda)$ (resp. $\Lambda=(0,0,\nu,0)$).  

It is easy to see that given $k_1 \in K_1$, we have 
$$k_1 \cdot (0,x,0,0)=(0,x,k_1x,0).$$ 
By writing  $(0,x,k_1x,0)=(0,x,0,0)(0,0,k_1x,0)(0,0,-\frac{k_1x^{2}}{2})$, we get that
$$[\rho_{\alpha,\lambda}(0,x,k_1x,0)f](u)=e^{-i \lambda k_1 xu+i \lambda k_1 \frac{x^{2}}{2}} f(u-x),$$
and hence we define
$$[\omega_{\alpha,\lambda}(k_1,0) f](u)=e^{-i \lambda \frac{u^{2}}{2}k_1} f(u).$$
Also, given $k_2 \in K_2$, 
$$k_2\cdot(0,x,0,0)=(0,x,0,k_2x).$$ 
Then 
$$[\rho_{\alpha,\lambda}(0,x,0,k_2x) f](u) = e^{i \lambda k_2 x} f(u-x),$$ 
and we set
$$[\omega_{\alpha,\lambda}(0,k_2)f](u)= e^{i \lambda k_2 u} f(u).$$
That is, 
$$[\omega_{\alpha,\lambda}(k_1,k_2)f](u)= e^{-i\lambda \frac{u^{2}}{2} k_1 + i \lambda k_2 u} f(u).$$\\

The analysis to the non generic orbit is similar, and we obtain
$$[\rho_{\nu}(0,x,k_1x,0)f](u) = e^{i \nu k_1 x} f(u-x),$$
$$[\rho_{\nu}(0,x,0,k_2 x)f](u) = f(u-x).$$

Then, 
$$[\omega_{\nu}(k_1,0)f](u) = e^{i \nu k_1 u} f(u),$$
$$[\omega_{\nu}(0,k_2)f](u) = f(u).$$
That is, 
$$[\omega_{\nu}(k_1,k_2)f](u) = e^{i \nu k_1 x} f(u).$$\\

It follows straightforward that \eqref{otros} holds for all $\Lambda \in \mathfrak{n}^{\ast}$ .\\


Thus we can conclude that the decomposition of $\omega_{\alpha, \lambda}$ on $L^{2}(\mathbb{R})$ is 
\begin{equation}
L^{2}(\mathbb{R})= \int_{\mathbb{R}} \chi_{-\lambda \frac{u^{2}}{2},\lambda u} \ du
\end{equation}
where $\chi_{-\lambda \frac{u^{2}}{2},\lambda u}(k_1,k_2)=e^{i(-k_1\lambda \frac{u^{2}}{2}+\lambda u k_2)}$.

Analogously, the decomposition of $\omega_{\nu}$ on $L^{2}(\mathbb{R})$ is 
\begin{equation}
L^{2}(\mathbb{R})= \int_{\mathbb{R}} \chi_{\nu u,0} \ du
\end{equation}
where $\chi_{\nu u,0}(k_1,k_2)=e^{ik_1\nu u}$.\\


In the last case, $\Lambda \equiv 0$ thus $\mathfrak{m}_{\Lambda}=\mathfrak{n}$  
and  $\chi_{\Lambda} \equiv 1$. 
Hence $\rho_{\Lambda}$ is the trivial representation. 
This case concludes the analysis of the metaplectic representation obtaining that $\omega_{\Lambda}$ is multiplicity free for all $\Lambda \in \mathfrak{n}^{\ast}$. 
Then, we obtain the following result.
\begin{teo}\label{Teo}
	Let $N=S \ltimes H_1$ and $K \subseteq Aut_0(N)$ defined in \eqref{K}. Then, the pair $(K,N)$  is a generalized Gelfand pair. 
\end{teo}

\end{document}